\documentclass[10pt]{article}

\usepackage[psamsfonts]{amssymb}
\usepackage{amsfonts,euscript,wrapfig}
\usepackage{amsmath,amsthm,graphicx}

\title{On Weak Separation Property for Affine Fractal Functions.}

\author{A.K.B.~Chand \and 
A.V.~Tetenov\footnote{Supported by Russian Foundation of Basic Research project
13-01-00513} }

\begin{document}

\maketitle

{\bf \large 1. Introduction.}  \ \   Weak  separation property (WSP),  which was  developed  since  90-ies  in  papers  of C.Bandt\cite{SSS7}, K.-S.Lau and S.-M.Ngai \cite{Lau} and M.Zerner\cite{Zer}    remains  one  of  the  main tools  of  analyzing  dimension problems. In recent  years it proved to be useful for the  study  of geometric  structure  of  self-similar  sets \cite{Atet1} and rigidity of self-similar  structures \cite{Atet2}.

In this  short  note we apply  this  notion to the   theory  of   affine fractal interpolation  funtions.

\bigskip

Standard  definition (see \cite{Barnsley1986}, \cite{Barnsley1988}, \cite{Nav})  of   affine  fractal  function   \ $f:[a,b]\to \mathbb R$ deals with a  partition $a=x_0<x_1<...<x_m=b$ of  the  interval $[a,b]$ and a system $\EuScript S=\{S_1,...,S_m\}$ of   affine  transformations
$$S_i(x,y)=\left(\begin{array}{cc}  p_i &0 \\ r_i &q_i \\  \end{array}\right) 
\left(\begin{array}{c} x \\ y \\  \end{array}\right) +
\left(\begin{array}{c}  h_i \\ s_i \\  \end{array}\right),   |p_i|<1,  |q_i|<1, $$
which  send  vertical  strip $a\le x\le b$ to vertical  strips $L_i=\{(x,y):x_{i-1}\le x \le x_i\}$. These  strips divide  the  graph $\Gamma(f)$ to {\bf non-overlapping} pieces $\Gamma_i=S_i(\Gamma(f))=\Gamma(f)\cup L_i$ whose  union  is $\Gamma(f)$.

\bigskip

But a more  general  approach  must  take  into account 
  the  possibility  of  overlaps:

For  example, a  system $\EuScript S$ consisting of  4 maps

$$ S_1:\left(\begin{array}{cc}  1/5 &0 \\  & \\ 1/5 &a \\  \end{array}\right) 
\left(\begin{array}{c} x \\ y \\  \end{array}\right),  \  \ \  
S_2:\left(\begin{array}{cc}  1/3&0 \\  & \\ -1/5&-1/5 \\  \end{array}\right) 
\left(\begin{array}{c} x \\ y \\  \end{array}\right) +
\left(\begin{array}{c} 1/5 \\   \\ 1/5\\  \end{array}\right),$$  $$
S_3:\left(\begin{array}{cc}   1/3 &0 \\  & \\  1/5&-1/5 \\  \end{array}\right) 
\left(\begin{array}{c} x \\ y \\  \end{array}\right) +
\left(\begin{array}{c}  {7}/{15}  \\    \\  0\\  \end{array}\right), \ \ 
S_4:\left(\begin{array}{cc}  1/5 &0 \\  & \\  -1/5&a \\  \end{array}\right) 
\left(\begin{array}{c} x \\ y \\  \end{array}\right) +
\left(\begin{array}{c}  {4}/{5}  \\   \\ 1/5\\  \end{array}\right)$$

\bigskip

defines \ \  a  self-affine  function  \ whose  graph \  passes  through \  the points\bigskip \\
$(0,0)$,\ $(1/5,1/5)$,\ $(7/15,0)$,\  $(8/15,0)$,\  $(4/5,1/5)$,\  $(1,0)$ and has  overlapping  pieces
$$S_2(\Gamma)\cap S_3(\Gamma)=S_2S_4(\Gamma)=S_3S_1(\Gamma)=\Gamma\left
(f|_{[\frac{7}{15} ,\frac{8}{15} ]}\right)$$

\begin{figure}[h] {\centering\includegraphics[scale=.33]{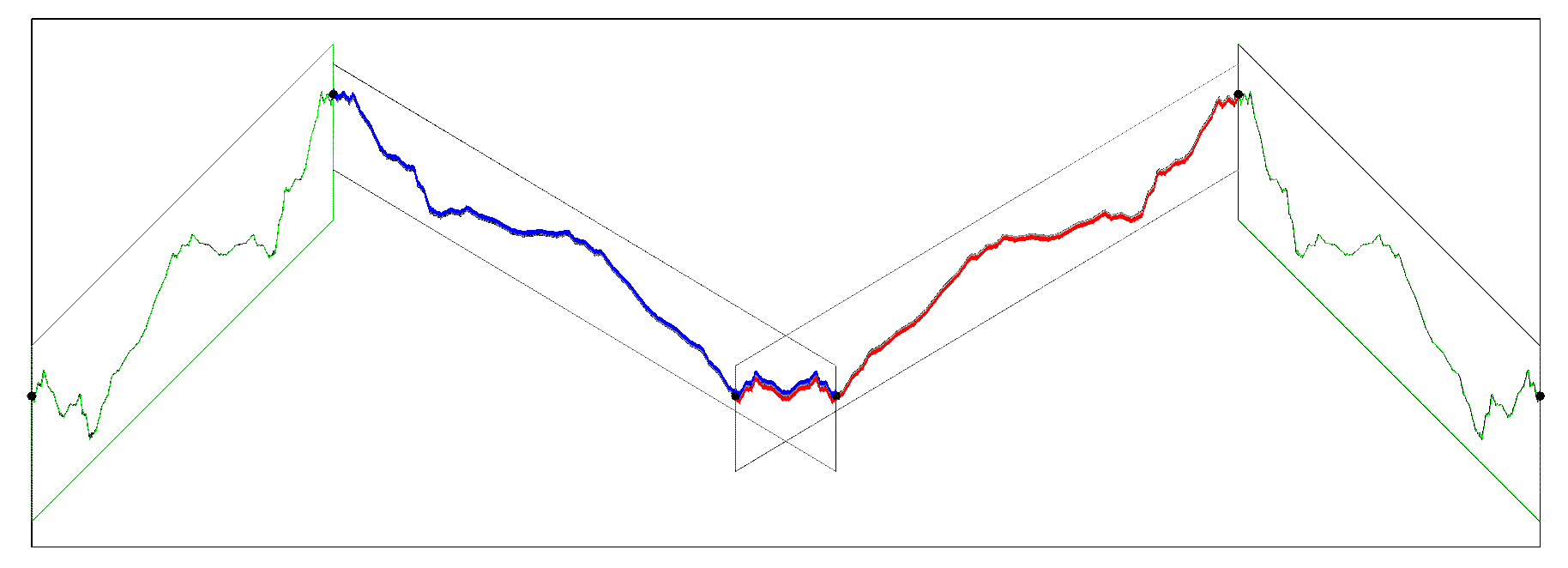}}
\caption{ The  graph of  $f(x)$: overlapping  pieces are    blue  and  red.}
\end{figure}

In view  of  the above  argument, we  use  the  following  definition which allows  the  overlaps:\\

\medskip

{\bf Definition 1.} {\it 
Let ${\EuScript S}=\{S_1, ... , S_m \}$ be a system of  affine  maps \\
$$S_i(x,y)=(p_ix+h_i, q_iy+r_ix+s_i); |p_i|,|q_i|<1$$\\
A function $f(x)$ is {\em affine  fractal  function} on  $[a,b]$ defined  by the  system  ${\EuScript S}$, if its graph $\Gamma(f)=\{(x,f(x)), x\in[a,b]\}$ is the  attractor  of   the  system  ${\EuScript S}$.} 

\bigskip

To formulate  the  main Theorem  we  recall some  definitions and  notation:

\bigskip

 We  denote  the projections  of $S_i$ to $\mathbb R$  by $S_i^{\diamond}(x)=p_ix+h_i$ and  we   denote 
 ${\EuScript S^\diamond}=\{S^\diamond_1, ... , S^\diamond_m \}$.\\
 $G$   denotes  the  semigroup generated  by ${\EuScript S}$ and ${G^\diamond}$ denotes  the  semigroup generated  by ${\EuScript S^\diamond}$.  

Observe that each  element  $g_{\bf i}=S_{i_1}S_{i_2}...S_{i_k}$ of the  semigroup  $G$ is  a  map of  the  form $g_{\bf i}(x,y)=(p_{\bf i}x+h_{\bf i}, q_{\bf i} y+r_{\bf i} y+s_{\bf i} ); |p_{\bf i} |,|q_{\bf i} |<1$  where $p_{\bf i}=p_{i_1}p_{i_2}...p_{i_k}$, $q_{\bf i}=q_{i_1}q_{i_2}...q_{i_k}$, while $g^\diamond_{\bf i}(x)=S^\diamond_{i_1}S^\diamond_{i_2}...S^\diamond_{i_k}(x)=(p_{\bf i}x+h_{\bf i})$.

\bigskip

We define {\it associated  families} $\EuScript F= G^{-1}\circ G$ and $\EuScript F^\diamond= G^{\diamond-1}\circ {G^\diamond}$ for  the  system ${\EuScript S}$ (resp. ${\EuScript S^\diamond})$. 
Each  element  of  the   family $\EuScript F$ is  a  composition  $g=g^{-1}_{\bf j}g_{\bf i}$ and also has  the  form $g(x,y)=(p x+h, q y+r x+s) $, while  its  projection  $g^\diamond=g^{\diamond-1}_{\bf j}g^\diamond_{\bf i}$  satisfies  $g^\diamond(x)=p x+h$.

\medskip

{\bf Definition 2.} { \it The  system ${\EuScript S}$ satisfies {\em weak separation property    (WSP)} if ${\rm Id}$ is  an isolated  point  in  the  associated  family $\EuScript F$.}

\medskip
So, if the  system $\EuScript S$ does  not  satisfy  the  weak separation property, there  is  a sequence $g_n\in \EuScript F$ which converges  to $\rm Id$.

\bigskip

{ \bf  \large 2. The  main Theorem.}

\bigskip

In this  paper  we  prove  the  following  Theorem:

\bigskip

{\bf Theorem 3.} { \it Let $f(x)$ be the affine  fractal  function defined by  a  system ${\EuScript S}$ on the  segment $[a, b]$.  If ${\EuScript S^\diamond}$ does  not  satisfy  weak  separation property, then the  graph $\Gamma(f)$ is  a segment of a parabola.}

\bigskip

First  of  all, it follows from the definition   that   each  fractal affine  function  is  continuous, because its  graph $\Gamma(f)$  is a compact set. \\

 Second,  a  remarkable property  of
 the  maps $g\in  {\EuScript F}$ is  that  these maps  move  the  points  of $\Gamma(f)$ along  $\Gamma(f)$:

\bigskip

{\bf Lemma 4.}  {\it  If 
 $g\in  {\EuScript F}$ and for some $x\in [a,b]$, $g^\diamond(x)\in [a,b]$, then $g(x,f(x))=(g^\diamond(x), f(g^\diamond(x))$.}

\bigskip

{\bf Proof.}    Let $g\in  {\EuScript F}$, so $g=g^{-1}_{\bf j}g_{\bf i}$. If $(x,y)\in \Gamma(f)$, then $ g_{\bf i}(x,y)\in \Gamma(f)$.  Suppose $(u,v)\in \Gamma(f)$ and
$g^\diamond_{\bf j}(u)=g^\diamond_{\bf i}(x)$. Since $g_{\bf j}(u,v)\in \Gamma(f)$, $g_{\bf j}(u,v)= g_{\bf i}(x,y)$, therefore $g^{-1}_{\bf j}g_{\bf i}(x,y)=(u,v)$, so $g(x,f(x))=(g^\diamond(x), f(g^\diamond(x))$. $\blacksquare$

\bigskip

These  facts  imply  that weak separation property   holds for  both  systems ${\EuScript S^\diamond}$  and ${\EuScript S}$ simultaneously:

\bigskip

{\bf Lemma 5.} { \it Let $f(x)$ be the affine fractal  function defined by  a  system ${\EuScript S}$ on the  segment $[a, b]$ whose  graph is  not a  straight line  segment.  The  system ${\EuScript S}$   satisfies  WSP iff  ${\EuScript S^\diamond}$   satisfies  WSP.}

\bigskip

{\bf Proof.}   
 Suppose   that WSP does  not  hold  for $ {\EuScript F^\diamond}$.
Take   three  points $(x_i,y_i),  i=1,2,3$ on $\Gamma(f)$ which do  not  lie  on a line. 
If $g_n^\diamond\to Id$ then for each $i$,  $g_n^\diamond(x_i)\to x_i$. Since $f$ is  continuous,
$g_n (x_i,y_i)\to (x_i,y_i)$. This means that $g_n$ converges to $Id$ and WSP does  not  hold  for $ {\EuScript F}$. 

Suppose now  that WSP does  not  hold  for $ {\EuScript F}$ and there  is a sequence $g_n\in  {\EuScript F}$ which  converges  to  $Id$.  Consider  the coefficients of $g_n(x,y)=(p_n x+h_n, q_n y+r_n  x+s_n)$:  $p_n$ and $q_n$ converge to 1,  while  $h_n, r_n$ and  $s_n$  converge to 0. Therefore $g^\diamond_n(x)=p_n x+h_n$ also converges to $Id$.     $\blacksquare$

\medskip

{\bf Lemma 6.} { \it Suppose $U$ is  a  family of  functions $\varphi(x)\in C^3[a,b]$, satisfying inequality $|\varphi(x)|\le M$. If for any $\varphi(x)\in U$, $\varphi''(x)$ and $\varphi'''(x)$ are monotonous and do not change  the  sign on $[a,b]$, then for  any segment $[a',b']\subset (a,b)$, the family $U'=\{\varphi|_{[a',b']}, \varphi\in U\}$ is bounded 
in $C^3([a',b'])$}

{\bf Proof.}  Without loss of generality,  we  suppose $[a,b]=[0,1]$ and $\varphi''(x)>0$ on $[0,1]$. Take some $\lambda\in(2^{1/3},1)$. \\

Since $\varphi(1)\le M$ and $\varphi(\lambda)\ge-M$, $\varphi'(\lambda)<\dfrac{2M}{1-\lambda}$. Similarly, we
get $\varphi'(1-\lambda)>\dfrac{-2M}{1-\lambda}$. So $\varphi'(x)<\left|\dfrac{2M}{1-\lambda}\right|$ on $[1-\lambda,\lambda]$.

\bigskip

Repeating the same step for $\varphi'$   we get $0<\varphi''(\lambda^2)<\dfrac{4M}{\lambda(1-\lambda)^2}$
 if  $\varphi''$ inreases  and
$\varphi''(1-\lambda^2)>\dfrac{4M}{\lambda(1-\lambda)^2}$ if  $\varphi''$ decreases, so 
$\varphi''(x)<\dfrac{4M}{\lambda(1-\lambda)^2}$ on $[1-\lambda^2,\lambda^2]$.

\bigskip

The  same  way,  we  have $|\varphi'''(x)|<\dfrac{8M}{\lambda^3(1-\lambda)^3}$ on $[1-\lambda^3,\lambda^3]$.
Taking such $\lambda$, that $[a',b']\subset [1-\lambda^3,\lambda^3]$, we  obtain  the  statement  for  the  segment 
$[a',b']$. $\blacksquare$  

\bigskip

{\bf Lemma 7.} { \it  Let $g\in  \EuScript F$ and  ${\rm fix}({g^\diamond})\notin [a,b]$.
Suppose that\\
(i)  \ \  if $x_1, x_2\in[a,b]$ and $|x_1-x_2|<\delta $,  then   $\|(x_1, f(x_1))-(x_2, f(x_2))\|<\varepsilon$;\\
(ii)\ \ \  $\|g(x,y)-(x,y)\|<\delta$ for  any  point $(x,y)\in \Gamma(f)$.\\ Then for  some $M\in \mathbb N$ either $\{g^n(a,f(a)), n=0,...,M\}$  or $\{g^n(b,f(b)),  n=0,...,M\}$   is  an $\varepsilon$-net  in $\Gamma(f)$.}

\bigskip

{\bf Proof.}   The  condition  (i) implies  that if $\{x_1,...,x_k\}$ is  a $\delta$-net  in $[a,b]$, then $\{ (x_1,f(x_1)),...,(x_k,f(x_k))\} $  is  an   $\varepsilon$-net in $\Gamma(f)$. So we  have to  show  that  $g^{\diamond n}(a)$ or  $g^{\diamond n}(b)$ form a $\delta$-net  in $[a,b]$.

Since ${\rm fix}({g^\diamond})\notin [a,b]$, we have either  $g^\diamond(x)>x$ for  any $x\in [a,b]$ or  $g^\diamond(x)<x$ for  any $x\in [a,b]$.

\medskip

Suppose $g^\diamond(a)>a$. Then for  any  point 
$x\in [a,b]$,  $g^\diamond(x)>x$ and  $g^\diamond(x)-x<\delta$. Since  the  limit point of  the  sequence  $g^{\diamond n}(a)$ is  outside $[a,b]$,  there  is   such $M\in \mathbb N$ for  which $g^{\diamond M}(a)<b<g^{\diamond M+1}(a)$, so for  any $n=1,...,M$, $g^{\diamond n}(a)-g^{\diamond n-1}(a)<\delta$ and $b-g^{\diamond M}(a)<\delta$. Therefore   $\{g^n(a,f(a)), n=0,1,....,M\}$ is  an $\varepsilon$-net in $\Gamma(f)$.  The proof  in  the  case $g^\diamond(b)<b$ is similar. $\blacksquare$

\bigskip

{\bf Lemma 8.} { \it Suppose $g(x,y)\in \EuScript F$,  $\rm{fix}({g^\diamond})\notin [a,b]$ and ${g^\diamond}(x)>x$ on $[a, b]$. Let $g^{\diamond T}(a)=b$. Then  the  set 
$\{g^t(a,f(a)), t\in[0,T]\}$ is  a  graph of  one of  the  following  functions  on $[a,b]$:\\
1. $y=Ax^2+Bx+C$;\\
2. $y=Ax+B e^{Kx}+C$;\\
3. $y=Ax+B(\log(x-C))+D$,  $C\notin[a,b]$.\\
4. $y=Ax+B(x-C)^K+D$, $C\notin[a,b]$;\\
5.  $y=A(x-C)\log(x-C) +Dx+E$,   $C\notin[a,b]$.
}\\

\bigskip

{\bf Proof}.   It  is  sufficient  to  check  the   statement  in the  case $a=0$, $b=1$, $f(0)=0$ and $p>1$. Since  $g$ is  close  to $Id$, $p$ and $q$ are  close  to 1 and  therefore  positive. 

The  five types  of  functions    arise    from  direct  solution  of  recurrence  equations :

\bigskip

 1.  If $g(x,y)=(x+h, y+r x+s)$, then  the  points $g^n(0,0)$ lie  on a  parabola  $y=Ax^2+Bx$, where $A=\dfrac{r}{2h}$ and $B=\dfrac{2s-hr}{2h}$;\\
\medskip

2.  If $g(x,y)=(x+h, qy+r x+s)$, $q\neq1$,  then  the  points $g^n(0,0)$ lie  on a  graph
of a function  $y=Ax+B (e^{Kx}-1)$, where $K=\dfrac{\log q}{h}$, $A=\dfrac{r}{q-1}$, 
$B=\dfrac{h r+(q-1) s}{(q-1)^2}$;

\medskip

3.  If $g(x,y)=(p x+h, y+r x+s)$, then the  points $g^n(0,0)$ lie  on a  graph of $y=Ax+B(\log(1+x/C))$, where
$C=\dfrac{h}{p-1}$, $A=\dfrac{r }{p-1}$,  $B=\dfrac{h r+(1-p) s}{(1-p)\log p}$;

\medskip

4.  If $g(x,y)=(p x+h,q y+r x+s)$, then the  points $g^n(0,0)$ lie  on a  graph of a function $y=Ax+B(x/C+1)^K-B$, where $A=\dfrac{r}{p-q}$, $C=\dfrac{h}{p-1}$, 
$B=\dfrac{h r+s (q-p)}{(q-1) (q-p)}$ and
$K=\dfrac{\log q}{\log p}$;

\medskip

5.  If $g(x,y)=(p x+h,p y+r x+s)$, then the  points $g^n(0,0)$ lie  on a  graph of a function $y=A(x/C+1)\log(x/C+1)+Bx$, where  $C=\dfrac{h}{p-1}$,  $A=\dfrac{rC}{p\log p}$
 and $B=\dfrac{C r-s}{C-C p}$.

\medskip

Applying to x coordinate a linear  transformation which  sends [a,b] to [0,1], we get  the  formulas 1-5 of  the  statement. $\blacksquare$

\bigskip

{\bf Proof of  the  Theorem 1.}
Take  such  sequence $g_n\to\rm Id$, $g_n\in \EuScript F$ and  such  segment 
$[a_1,b_1]\subset (a,b)$, that  for  any $n$, ${\rm fix}(g^\diamond_n)\notin[a_1,b_1]$.

Since $g^{-1}_n$ also  converge  to $\rm Id$, we  may  suppose that for  any $n$,  $p_n\ge 1$. 

Without  loss  of  generality  we may  suppose that for  any $n$, $g^\diamond_n(a_1)>a_1$. 
Let $T_n$ be  such  number, that $g_n^{\diamond T_n}(a_1)=b_1$. 
Each  curve $\{g_n^t(a_1,f(a_1)), t\in [0,T_n]\}$ is  a graph  of  a  function $\varphi_n(x)$ on  the  segment $[a_1,b_1]$.

It  follows  from  Lemma 7 that $\varphi_n(x)$ uniformly  converges  to $f(x)$ on $[a_1,b_1]$.

By Lemma 8, each  of  these  functions is of  one  of  5 types, indicated  by  the Lemma.
Therefore  the  functions $\varphi_n(x)$ have  monotonous derivatives  $\varphi''_n(x)$, $\varphi'''_n(x)$, which do  not  change  their  sign on $[a_1,b_1]$. By Lemma 6, for  any
$[a_2,b_2]\in(a_1,b_1)$,
 the  family $\{\varphi_n(x)|_{[a_2,b_2]}\}$ is a  bounded  subset  of $C^3([a_2,b_2])$.
Therefore  some  subsequence of $\varphi_n(x)$ converges  in $C^2([a_2,b_2])$, which implies  that  $f(x)$ is  twice differentiable  on $[a_2,b_2]$.  This  means  that $f(x)\in C^2([a,b])$. As it  was proved  in \cite[Theorem 3]{BKra} this  implies that $\Gamma(f)$ is  a  parabolic segment. $\blacksquare$


\bigskip



 
\begin{thebibliography}{11}

\bibitem{BKra} Ch.~Bandt, A.~Kravchenko. Differentiability of fractal curves,  Nonlinearity, 2011, 24, pp.~2717- 2728.

\bibitem{SSS7}
   Ch.~Bandt, S.~Graf.  Self-similar sets 7. A characterization of
 self-similar fractals with positive Hausdorff measure. // Proc. Amer.
 Math. Soc., 1992, V.~114, No.~4, pp.~995--1001.

\bibitem{Barnsley1986}
{M.F.~Barnsley.}, Fractal functions and interpolation // Constructive Approximation, 1986,  V.~2, pp.~303-329.

\bibitem{Barnsley1988}
{M.F.~Barnsley.} Fractals Everywhere // Academic Press, 1988, 394 p. 

\bibitem{Hut}
J.~Hutchinson. Fractals and self-similarity. // Indiana Univ.
Math. J., 1981, V.~30, No.~5,  pp.~713--747.

\bibitem{Lau}
K.-S. Lau and S.-M. Ngai, Multifractal measures and a weak separation condition,// Adv. Math., 1999, No.~141, pp.~45-96.

\bibitem{Nav}
Navascu\'es,~ M.A., Chand,~ A.K.B., Veedu,~V.P. and Sebasti\'an,~M.V. // Fractal Interpolation Functions: A Short Survey. Applied Mathematics, (2014), 5, 1834-1841.

\bibitem{Atet1}  A.~V.~Tetenov. Self-similar Jordan arcs and
 graph-directed systems of similarities, Siberian.Math.J., 47, No.~5 (2006),
 pp.~940-949.

\bibitem{Atet2}  A.~V.~Tetenov,  On the rigidity of one-dimensional systems of
contraction similitudes., Siberian Electronical Math. Reports,V.3
(2006),  pp.~342--345.

\bibitem{Zer}   M.~P.~W.~Zerner, Weak separation properties for self-similar
sets.// Proc.Amer.Math.Soc. 1996,  V. 124, No.~11, pp.~3529--3539.



\end{thebibliography}
\end{document}